\documentclass[11pt,a4paper]{article}

\title{On shelling and flag vectors}

\author{Jonathan Fine\relax
\thanks{203 Coldhams Lane, Cambridge, CB1 3HY, England.
\quad E-mail: \texttt{j.fine@pmms.cam.ac.uk}}%
}
\date{1 October 1997}

\newcommand\ftilde{\widetilde{f}}
\newcommand\bibrule{\rule{2pc}{0.4pt}}

\begin{document}
\maketitle

\begin{abstract}\noindent
This note defines a flag vector for $i$-graphs.  The construction applies
to any finite combinatorial object that can be shelled.  Two possible
connections to quantum topology are mentioned.  Further details appear in
the author's \emph{On quantum topology, hypergraphs and flag vectors},
(preprint q-alg/9708001).
\end{abstract}

\noindent
The purpose of this note is to state concisely a new combinatorial
concept.  This note is an announcement; details appear elsewhere
\cite{bib.JF.QTHGFV}.  Throughout, suppose that $G$ is something that
(a)~is built out of cells, and (b)~can be shelled.  To fix ideas, it will
be assumed that $G$ is an $i$-graph on a finite vertex set $V$.  In other
words, $G$ is a possibly empty collection of $i$-element subsets of $V$. 
These subsets are the cells, also known as edges.

When a vertex $v$ is removed from $V$, the cells that contain $v$ must be
removed from $G$.  Each of these cells has exactly $i$ elements, one of
which is $v$ itself.  Thus, the link $L_v$ of $G$ at $v$ is defined to be
the $(i-1)$-graph, whose cells are exactly the $(i-1)$-subsets of $V -v$
which, upon the addition of $v$, become cells of $G$.

To shell $G$ is to remove the vertices $v$ from $V$ one at a time, until
none are left.  As this is done, a record is kept of the resulting links
$L_1$, $L_2$, $\dots$, $L_{N-1}$, $L_N$.  Here, $N$ is the number
of elements in the vertex set $V$.  Each shelling thus determines a
sequence $\{L_j\}$ of $(i-1)$-graphs.

The shelling vector $\ftilde G$ is defined inductively as follows.  Assume
$\ftilde L_i$ is defined, for $(i-1)$-graphs.  The sum over all shellings
\[
    \ftilde G = \sum \nolimits _{\rm shellings} \>
        \ftilde L _ 1 \otimes \ldots \otimes \ftilde L _ N
\]
is the shelling vector of $G$.  The induction is founded on $0$-graphs. 
There is only one zero element set, namely the empty set, and it is a
subset of any vertex set.  Let $a=a_V$ and $b=b_V$ be the two possible
$0$-graphs on the vertex set $V$.  The $0$-graph $a$ has no cells, while
$b$ has the empty set as a cell.  Use the equations $\ftilde a = a$,
$\ftilde b = b$ to found the induction.

The shelling vector is for various reasons too large.  The flag vector
$fG$ is defined inductively as follows.  Much as before, one writes
\[
    fG = \sum \nolimits _{\rm shellings} \>
        f' L _ 1 \otimes \ldots \otimes f' L _ N
\]
where now each $f'L_j$ is a linear function of the flag vector of each
link appearing in the shelling.  The induction is founded on $0$-graphs
exactly as before.

The link contributions $f'L$ are defined as follows.  Let $A$ and $B$ be
two cells of $L$, that do not have a common vertex.  Let $L_{{+}{+}}=L$,
$L_{{-}{+}}$, $L_{{+}{-}}$, $L_{{-}{-}}$ be the four graphs obtained by
removing none, one or both of $A$ and/or $B$ respectively from $L$.  The
flag vector $f$ takes $L$ to a point $fL$ lying in some vector space $W$. 
Define $W'$ to be $W$ modulo the span of all expressions of the form
\[
    fL_{{+}{+}} \> - \> fL_{{-}{+}} \> - \> 
    fL_{{+}{-}} \> + \> fL_{{-}{-}}
\]
and then set $f'L$ to be the residue of $fL$ in $W'$.  This completes the
definition of $i$-graph flag vectors.

In the same way, a flag vector can be defined for any finite object, built
out of cells, that can be shelled.  For ordinary or $2$-graphs, the flag
vector of graphs on $N$ vertices has $p(N)$, the number of partitions of
$N$, independent components.  It satisfies subtle, and as yet unexplored,
inequalities.  The author hopes that techniques similar to those described
in \cite{bib.JF.CPLA} will make these inequalities accessible. Each
triangulated $n$-manifold determines an $(n+1)$-graph.  It is hoped that
the quantum invariants of a manifold can be expressed as linear functions
of the flag vector (of any of its triangulations).  All this is discussed
further in \cite{bib.JF.QTHGFV}. The `disjoint pair of optional cells'
rule used to define $f'L$ from $fL$ is similar to the `independent
regional change' concept in Vassiliev theory~\cite{bib.JF.VTRCKI}.


\begin{thebibliography}{1}

\bibitem{bib.JF.QTHGFV}
J. Fine, Quantum topology, hypergraphs and flag vectors,
preprint q-alg/9708001 (August 1997)

\bibitem{bib.JF.CPLA}
\bibrule, Convex polytopes and linear algebra,
preprint alg-geom/9710001 (October 1997)

\bibitem{bib.JF.VTRCKI}
\bibrule, Vassiliev theory, regional change, and the Kontsevich integral,
(in preparation)

\end{thebibliography}
\end{document}